\documentclass[12pt,leqno]{amsart}
\theoremstyle{plain}
 \newtheorem{theorem}{Theorem}
 \newtheorem*{theorem*}{Main Theorem}
 
 \newtheorem{lemma}[theorem]{Lemma}
 \newtheorem*{lemma*}{Main Lemma}

\theoremstyle{definition}
 \newtheorem{definition}[theorem]{Definition}

 \newtheorem*{example*}{Example}
\newcommand\ord{\operatorname{ord}}
\newcommand\pd{\partial}
\newcommand\la{\langle}
\newcommand\ra{\rangle}
\newcommand\C{\mathbb{C}}
\newcommand\R{\mathbb{R}}
\newcommand\Z{\mathbb{Z}}
\newcommand{\Q}{\boldsymbol{Q}}
\newcommand{\CP}{\mathbb{CP}}
\newcommand\abold{\boldsymbol{a}}
\newcommand\ebold{\boldsymbol{e}}
\renewcommand\Re{\operatorname{Re}}
\renewcommand\Im{\operatorname{Im}}
\newcommand{\TC}{\operatorname{TC}}

\makeatother
\title{
   Minimal surfaces that attain equality\\
   in the Chern-Osserman inequality
}
\author{M. Kokubu}
\address[Masatoshi Kokubu]{%
   Department of Natural Science,
   Tokyo Denki University,
   Inzai, Chiba 270-1382, Japan%
}
\email{kokubu@chiba.dendai.ac.jp}
\author{M. Umehara}
\address[Masaaki Umehara]{%
   Department of Mathematics, Faculty of Science,
   Hiroshima University,
   Higashi-Hiroshima 739-8526, Japan%
}
\email{umehara@math.sci.hiroshima-u.ac.jp}
\author{K. Yamada}
\address[Kotaro Yamada]{%
   Faculty of Mathematics,
   Kyushu University 36, 
   Fukuoka 812-8185, Japan%
}
\email{kotaro@math.kyushu-u.ac.jp}
\subjclass{Primary 53A10; Secondary 53A07, 53C42}
\begin{document}
\begin{abstract}
In the previous paper, Takahasi and the authors generalized the theory
 of minimal 
surfaces in Euclidean $n$-space to that of surfaces with holomorphic
 Gauss map in
certain class of non-compact symmetric spaces. 
It also includes the theory of constant
mean curvature one surfaces in
 hyperbolic 3-space. Moreover, a Chern-Osserman type
inequality for such surfaces was shown. Though its equality condition
 is not solved yet, 
the authors have noticed that the equality condition of the original
 Chern-Osserman
inequality itself is not found in any
 literature except for the case $n=3$, in spite
of its importance.
In this paper, a simple geometric  
condition for minimal surfaces that attains
equality in the Chern-Osserman inequality is given. 
The authors hope it will be a useful reference for readers.
\end{abstract}
\maketitle
The total curvature $\TC(M)$ of any complete minimal surface $M$ in $\R^n$ 
has a value in $2\pi \Z$ and satisfies the following inequality 
called the {\it Chern-Osserman inequality\/} \cite{c-o}:
 \begin{equation}\label{eq:ch-os}
  \TC(M) \le 2\pi(\chi_M-m), 
 \end{equation}
where $\chi_M$ denotes the Euler number of $M$ and $m$ is the 
number of ends of $M$. 

Then it is natural to ask which surfaces attain the equality of the 
inequality \eqref{eq:ch-os}. 
In the case of $n=3$, Jorge and Meeks \cite{jm} gave a geometric 
proof of \eqref{eq:ch-os} and proved that the equality holds if 
and only if all of the ends are embedded.
However, for general $n>3$, the authors do not know any references
on it. 
The purpose of this paper is to give the following geometric 
condition for attaining equality in \eqref{eq:ch-os} for general $n$. 
\begin{theorem*}
 A complete minimal surface in $\R^n$ attains equality in the 
 Chern-Osserman inequality if and only if each end is asymptotic to 
 a catenoid-type end or a planar end in some $3$-dimensional subspace
 $\R^3$ in $\R^n$. 
 In particular,  all ends are embedded. 
\end{theorem*}
For $n=3$, according to Jorge-Meeks \cite{jm} and Schoen \cite{schoen},
one can easily observe that embedded ends are all
asymptotic to catenoids or planes (see Appendix).
So our theorem generalizes the result in Jorge-Meeks.
For general $n$ $(>3)$,
we remark that the embeddedness of ends is not a sufficient condition 
for the equality of \eqref{eq:ch-os}. 
For example, an embedded holomorphic curve 
$f \colon \C \setminus \{ 0 \} \to \C^2$ 
defined by ${f}(z)=(z,1/z^2)$  
(considered as a complete minimal surface in $\R^4$) 
has total curvature $-6\pi$. 
So it does not satisfy equality in \eqref{eq:ch-os}. 
%%%%%%%%%%%%%%%%%%%%%%%%%%%%%%%%%%%%%%%%%%%%%%%%%%%%%%%%%%%%%%%%%%%%%%
\section*{Preliminaries}
We shall review the properties of minimal surfaces in
$\R^n$ (cf. \cite{lawson}).
Let $f=(f_1,\dots,f_n)\colon M \to \R^n$ be a conformal minimal
immersion of a Riemann surface $M$, where $n\ge 3$ is an integer. 
Then $\pd f$ is a $\C^n$-valued holomorphic $1$-form on $M$.
We define the Gauss map $\nu\colon{}M\to \CP^{n-1}$ of $f$ as
\[
   \nu:=[\pd f]=\left[
                    \frac{\pd f_1}{\pd z}\colon{}
                    \frac{\pd f_2}{\pd z}\colon{}\cdots
                    \colon{}
                    \frac{\pd f_n}{\pd z}
                  \right],
\]
where $z$ is a complex coordinate of $M$.
Since $f$ is conformal, we have
\begin{equation}\label{eq:nullity}
   \la \pd f, \pd f \ra
       = \sum_{j=1}^n\left(\frac{\pd f_j}{\pd z}\right)^2\,dz^2=0.
\end{equation}
Thus, the Gauss map $\nu$  is valued in the complex quadric
$\Q^{n-2}\subset \CP^{n-1}$.

We assume that $f$ is complete and of finite total curvature. 
Under this assumption, the following properties are well-known: 
\begin{itemize}
\item $M$ is biholomorphic to a 
      compact Riemann surface $\overline M$ punctured at finitely many points 
      $\{p_1,\dots,p_m \}$. 
      Each point $p_j$ is called an {\it end}. 
\item The Gauss map $\nu$ can be extended holomorphically on $\overline M$, and 
      the total curvature is given by $-2\pi d$ where $d$ is the 
      homology degree of $\nu(\overline M)$ in $\CP^{n-1}$. 
\item For each end $p_j$, there exists 
      a local complex coordinate $z$ on $\overline M$ centered at $p_j$
      such that the first fundamental form $ds^2$ is 
      written as 
\[
      ds^2=|z|^{2\mu_j}dz\, d\bar z\qquad (\mu_j\le -2).
\]
      We call $\mu_j$ the {\it order\/} of the metric $ds^2$ at the
      end $p_j$ and denote by $\ord_{p_j}ds^2=\mu_j$.
      Since $ds^2=2\la\pd f,\overline{\pd}f\ra$,
      $\mu_j$ coincides with the order of $\pd f$ at the end $p_j$.
\end{itemize}
\begin{definition}\label{def:endtype} 
 An end $p_j$ of 
 $f\colon M={\overline M} \setminus \{p_1,\dots,p_m \} \to \R^n$ 
 is said to be asymptotic to a {\it catenoid-type\/} 
 (resp.~{\it planar\/}) end if there exists a piece of the catenoid 
 (resp.~the plane) 
\[
   f_0 \colon \{ |z-p_j| < \varepsilon \} \to \R^3 \subset \R^n 
\]
 which is complete at $p_j$ such that 
 $|f(z)-f_0(z)|=O(|z-p_j|)$, that is,
\[
 \frac{|f(z)-f_0(z)|}{|z-p_j|} 
\]
 is bounded on
 $\{|z-p_j| < \varepsilon \}$ for sufficiently small 
 $\varepsilon>0$.
\end{definition}

%%%%%%%%%%%%%%%%%%%%%%%%%%%%%%%%%%%%%%%%%%%%%
\section*{Proof of the Main Theorem}
The Chern-Osserman inequality
follows from the fact $\ord_{p_j}ds^2\le -2$ at each end $p_j$.
Moreover, equality holds if and only if $\ord_{p_j}ds^2=-2$ 
(see \cite[pp.~135--136]{lawson}, for example).
Thus the Main Theorem immediately follows from the following Lemma.
\begin{lemma}\label{lem:loc-eq-cond}
 Let $f\colon \Delta^* \to \R^n$ be a conformal minimal immersion 
 of a punctured disc 
 $\Delta^* = \{ z \in \C \ | \ 0 < |z| <1 \}$ into $\R^n$
 which is complete at the origin $0$. 
 Then $\ord_{0}ds^2=-2$ holds if and only if  the end $0$  is asymptotic
 to a catenoid-type end or a planar end in $\R^3 (\subset \R^n)$. 
 In particular, it is an embedded end. 
\end{lemma}
\begin{proof}
 Suppose that $\ord_{0}ds^2=-2$.   
 It implies that 
 the Laurent expansion of $\pd f$ is given by 
\begin{equation}\label{eq:laur1}
  \pd f=\left( \frac{1}{z^2}\abold_{-2}+\frac{1}{z}\abold_{-1} 
               + \cdots \right) dz, 
  \quad \abold_{-2} \in \C^n\setminus\{0\},\ \abold_{-1} \in \R^n
\end{equation}
 because the residue of $\pd f$ must be real.
 Moreover, it follows from \eqref{eq:nullity} that 
\[
 \la \abold_{-2}, \abold_{-2} \ra = 0, 
       \quad \text{and}\quad
 \la \abold_{-2}, \abold_{-1} \ra = 0. 
\]
 Therefore we have
\begin{gather*}
  |\Re \abold_{-2}|=|\Im \abold_{-2}|, \qquad 
  \la \Re \abold_{-2}, \Im \abold_{-2} \ra =0, \\
  \la \Re \abold_{-2}, \abold_{-1} \ra =0, \qquad 
  \la \Im \abold_{-2}, \abold_{-1} \ra =0. 
\end{gather*}
 Hence we can choose an orthonormal basis $\ebold_1,\dots,\ebold_n$ 
 of $\R^n$ so that 
\[
  \Re \abold_{-2} = a\, \ebold_1, \quad 
  \Im \abold_{-2} = a\, \ebold_2, \quad 
  \abold_{-1} = b\, \ebold_3 
\]
 for some real constants $a(\ne 0)$, $b$. 
 With respect to this basis, we have 
\[
  \pd f=\left( \frac{a}{z^2}(\ebold_1 + i \ebold_2) + 
  \frac{b}{z}\ebold_{3}+ \cdots \right) dz, 
  \quad a, b \in \R, (a \ne 0)
\]
 Then using the polar coordinate $z=r e^{i\theta}$, we have
\begin{equation}\label{eq:intoflaur1}
  f(z) = 2 \int_{z_0}^z \partial f 
       = -\frac{2 a\cos \theta}{r}\ebold_1-\frac{2 a\sin \theta}{r}\ebold_2 
         +2 b \log r \ebold_3 + O(r),
\end{equation}
 where $z_0$ is a base point.
 Here, we have dropped the constant terms in $f(z)$ by a suitable
 parallel translation. 
 By Definition~\ref{def:endtype}, 
 the formula \eqref{eq:intoflaur1} implies that the surface $f(\Delta^*)$ 
 is asymptotic to the catenoid (resp.~the plane) for the sufficiently 
 small $r$ if $b \ne 0$ (resp.~if $b=0$). 

 Conversely, suppose that $\ord_{0}ds^2 \ne -2$. 
 It implies that $\ord_{0}ds^2=-k$ ($k \ge 3$) and 
\begin{equation}\label{eq:laur3}
  \pd f=\left( \frac{1}{z^k}\abold_{-k}+ \cdots + \frac{1}{z}\abold_{-1} 
               + \cdots \right) dz, 
  \quad \abold_{-k} \ne 0 \in \C^n,\ \abold_{-1} \in \R^n. 
\end{equation}
 It is obvious that the end is asymptotic to neither a catenoid-type end 
 nor a planar end. 

 From now on, we shall prove that an end is embedded  
 if it is asymptotic to a catenoid-type end or a planar end. 
 Assume that the end is not embedded. 
 Then there exist two sequences $\{z^{}_j\}$, $\{z'_j\}$ convergent to 
 $0$ such that $f(z^{}_j)=f(z'_j)$ for all $j$. 
 Then by \eqref{eq:intoflaur1}, 
 there exists a positive constant $C$ such that 
\[
   \left| \frac{\cos \theta_j}{r^{}_j}-\frac{\cos \theta'_j}{r'_j} \right| 
         \le C |r^{}_j - r'_j |, \qquad 
   \left| \frac{\sin \theta_j}{r^{}_j}-\frac{\sin \theta'_j}{r'_j} \right| 
      \le C |r^{}_j - r'_j |,
\]
 where $z^{}_j = r^{}_j e^{i\theta_j}$ and $z'_j=r'_j e^{i\theta'_j}$
 ($j=1,2,\dots$).
 With these estimates, we have
\begin{align*}
 \left( \frac{1}{r^{}_j} - \frac{1}{r'_j} \right)^2 & \le 
 \frac{1}{{r^{}_j}^2} + \frac{1}{{r'_j}^2}
 -\frac{2}{r^{}_jr'_j}\cos(\theta_j-\theta'_j) \\ 
 & = \left| \frac{\cos \theta_j}{r^{}_j}-\frac{\cos \theta'_j}{r'_j} \right|^2 
 + \left| \frac{\sin \theta_j}{r^{}_j}-\frac{\sin \theta'_j}{r'_j} \right|^2 \\ 
 & \leq  2 C^2 |r^{}_j - r'_j |^2 ,
\end{align*}
 and then, 
\begin{equation}\label{eq:rnrn}
 \frac{1}{(r^{}_jr'_j)^2} \le 2C^2
\end{equation}
 holds.
 However the left hand side of \eqref{eq:rnrn} diverges to $+\infty$ as 
 $j \to \infty$. This is a contradiction. 
\end{proof}
Besides the Chern-Osserman inequality \eqref{eq:ch-os}, 
the following inequalities for fully immersed 
complete minimal surfaces are known. 
(We say that the immersion $f$ is {\it full\/} if 
 the image $f(M)$ is not contained in any hyperplanes of $\R^n$.)

Gackstatter \cite{gack} proved that 
\[
   \TC(M) \le (2\chi_M+m-1-n)\pi.
\]

On the other hand, Ejiri \cite{ejiri} proved the inequality 
\begin{equation}\label{eq:ejineq}
    \TC(M) \le (\chi_M + m -2n + 2l)\pi
\end{equation}
if its Gauss image $\nu (M)$ is contained in an ($n-1-l$)-dimensional
subspace of $\CP^{n-1}$. 

Here, we shall give a new example of complete minimal surfaces 
which satisfies the equality both in the Chern-Osserman equality 
\eqref{eq:ch-os} and 
in the Ejiri inequality \eqref{eq:ejineq}. 
\begin{example*}[Generalized Jorge-Meeks' surface]
 For $j=0,1,\dots,m-1$, we put 
 \[
  g_{j}(z)=\frac{z^j(1-z^{2m-2j})}{(z^{m+1}-1)^2}, \quad 
   h_{j}(z)=\frac{iz^j(1+z^{2m-2j})}{(z^{m+1}-1)^2}, 
 \]
 and define a complete conformal minimal immersion by 
 \begin{equation}\label{eq:imm-gjm}
  f_m:=\Re \int_{z_0}^z\left(
                 g_{0},h_{0},g_{1},h_{1}, \ldots , g_{m-1},h_{m-1}, 
                 \frac{2 \sqrt{m} z^m}{(z^{m+1}-1)^2}
               \right)\,dz.  
 \end{equation}
 Then by similar computations as in \cite{jm}, the integrand of
 \eqref{eq:imm-gjm} has real residue at each pole, and then,
 $f_m$ gives a conformal minimal immersion
 \[
   f_m\colon{}M = \bigl(\C\cup\{\infty\}\bigr)\setminus
      \{ z\,;\,z^{m+1}=1\}\longrightarrow
      \R^{2m+1}.
 \]
 Obviously, the genus of $M$ is zero, the number of ends is $m+1$, 
 and $f_m \colon M \to \R^{2m+1}$ is full. 

 Since the degree of the Gauss map of $f_m$ is $2m$, 
 the total curvature  $\TC(M)$ is equal to $-4 m \pi$. 
 Therefore it attains the equality in the Chern-Osserman inequality. 
 
 On the other hand, it is easy to see that $f_m$ has non-degenerate 
 Gauss map, that is, $l=0$ in \eqref{eq:ejineq}. 
 Then the right hand side of \eqref{eq:ejineq} is $-4 m \pi$. 
 Hence the equality in \eqref{eq:ejineq} holds. 
\end{example*}
%%%%%%%%%%%%%%%%%%%%%%%%%%%%%%%%%%%%%%%%%%%%%%%%%%%%%%%%
\section*{Appendix: Embedded ends in $\R^3$}
 For the case $n=3$, embeddedness of the end $0$ in
 Lemma~\ref{lem:loc-eq-cond} implies $\ord_0ds^2=-2$, 
 and consequently the end is asymptotic to a catenoid-type end or
 a planer end
 (\cite[Theorem~4]{jm} or \cite[Proposition~1]{schoen}).
% We give a simple proof of it.
Here we shall give a simple proof of this fact, which is a mixture of 
Jorge-Meeks' and Schoen's. The authors hope that it will be 
 helpful to readers.
The crucial point of the Jorge-Meeks' proof is to show that the intersection
of the end and the sphere of radius $r$ centered at the origin
converges to a finite covering of a great sphere as $r\to \infty$.
According to Schoen \cite{schoen}, 
we prove it via the Weierstrass representation directly.

 Consider the Laurent expansion as \eqref{eq:laur3} for $k\geq 2$.
 Without loss of generality, we may set $\abold_{-k}=(a,ia,0)$
 $(a\in\R\setminus\{0\})$ because of \eqref{eq:nullity}.
 Integrating this, we have
\[
    f(re^{i\theta}) = \frac{1}{r^{k-1}}
             \left[
               2a \bigl(
                 \cos(k-1)\theta, \sin (k-1)\theta, 0
               \bigr)
                 + o(1)\right],
\]
 where $o(1)$ means a term tending to $0$ as $r\to 0$.
 Let $S^2_R$ be the sphere in $\R^3$ with radius $R$ centered at the
 origin and consider the intersection of the surface and $S^2_R$:
\[
   E_R:=\frac{1}{R}\left(
                     S^2_R\cap f(\Delta^*)
                   \right)\subset S^2_1,
\]
 which is normalized as a subset of the unit sphere.

 Here, $f\in S^2_R$ if and only if 
\[
   R^2=f_1^2 + f_2^2 +f_3^2 = \frac{1}{r^{2k-2}}(4a^2+o(1))
\]
 holds.
 Then $r\to 0$ as $R\to\infty$ when $f(r e^{i\theta})\in S^2_R$ because
 $k\geq 2$. In particular,
 $\lim_{r\to\infty}R^2 r^{2k-2}=4a^2$ 
 holds.

 Then under the condition $f(z)\in S^2_R$, 
\[
   \lim_{R\to\infty}\frac{1}{R}f(re^{i\theta})
       = \bigl(\cos(k-1)\theta,\sin(k-1)\theta,0\bigr)
\]
 holds.
 This implies that, for sufficiently large $R$,
 $E_R$ is a closed curve in a neighborhood of the 
 equator of $S^2_1$ with rotation index $|k-1|$, which is embedded
 if and only if $k=2$.\qquad \qedsymbol

\medskip
\noindent
{\bf Acknowledgement.}
We would like to thank Wayne Rossman for valuable comments.

%%%%%%%%%%%%%%%%%%%%%%%%%%%


\begin{thebibliography}{KTUY}
\bibitem[CO]{c-o}
S. Chern and R. Osserman, 
\newblock{\it Complete minimal surface
            in Euclidean $n$-space,}
\newblock{J. Analyse Math., {\bf 19} (1967) 15--34.}
%%
\bibitem[E]{ejiri}N.~Ejiri, 
\newblock{\it Degenerate minimal surfaces of finite total curvature in
            $R^N$,} 
\newblock{Kobe~J.~Math., {\bf 14} (1997), 11--22.}
%%%
\bibitem[G]{gack}F.~Gackstatter, 
\newblock{\em \"Uber die Dimension einer Minimalfl\"ache und zur Ungleichung 
            von St.~Cohn-Vossen,} 
\newblock{ Arch.~Rational Mech. Anal., {\bf 61} (1976), 141--152.}
%%
\bibitem[JM]{jm}L.~P.~M.~Jorge and W.~H.~Meeks~I\/I\/I, 
\newblock{\em The topology of complete minimal surfaces 
            of finite total curvature, }
\newblock{Topology, {\bf 22} (1983), 203--221.}
%%
\bibitem[KTUY]{ktuy} M.~Kokubu, M.~Takahashi, M.~Umehara and
K.~Yamada, 
\newblock{\em 
      An analogue of minimal surface theory
       in \boldmath$\mbox{SL}(n,\C)/\mbox{SU}(n)$, }
\newblock{Preprint.}
%
\bibitem[L]{lawson}
H.~B.~Lawson,
\newblock{\sc Lectures on minimal submanifolds (Volume 1),}
\newblock{Publish or Perish Inc., 1980.}
%%
\bibitem[S]{schoen}
R. Schoen, 
\newblock{\em Uniqueness, symmetry and embeddedness of minimal
       surfaces,}
\newblock{J. Differential Geometry, {\bf 18} (1983), 791--809.}
\end{thebibliography}
\end{document}